\documentclass[12pt]{article}
\usepackage{amsmath,amsfonts,amssymb,amsthm}
\input amssym.def
\begin{document}
\def \Z{\Bbb Z}
\def \C{\Bbb C}
\def \R{\Bbb R}
\def \Q{\Bbb Q}
\def \N{\Bbb N}

\def \A{{\mathcal{A}}}
\def \D{{\mathcal{D}}}
\def \E{{\mathcal{E}}}
\def \E{{\mathcal{E}}}
\def \H{\mathcal{H}}
\def \S{{\mathcal{S}}}
\def \wt{{\rm wt}}
\def \tr{{\rm tr}}
\def \span{{\rm span}}
\def \Res{{\rm Res}}
\def \Der{{\rm Der}}
\def \End{{\rm End}}
\def \Ind {{\rm Ind}}
\def \Irr {{\rm Irr}}
\def \Aut{{\rm Aut}}
\def \GL{{\rm GL}}
\def \Hom{{\rm Hom}}
\def \mod{{\rm mod}}
\def \ann{{\rm Ann}}
\def \ad{{\rm ad}}
\def \rank{{\rm rank}\;}
\def \<{\langle}
\def \>{\rangle}

\def \g{{\frak{g}}}
\def \h{{\hbar}}
\def \k{{\frak{k}}}
\def \sl{{\frak{sl}}}
\def \gl{{\frak{gl}}}

\def \be{\begin{equation}\label}
\def \ee{\end{equation}}
\def \bex{\begin{example}\label}
\def \eex{\end{example}}
\def \bl{\begin{lem}\label}
\def \el{\end{lem}}
\def \bt{\begin{thm}\label}
\def \et{\end{thm}}
\def \bp{\begin{prop}\label}
\def \ep{\end{prop}}
\def \br{\begin{rem}\label}
\def \er{\end{rem}}
\def \bc{\begin{coro}\label}
\def \ec{\end{coro}}
\def \bd{\begin{de}\label}
\def \ed{\end{de}}

\newcommand{\m}{\bf m}
\newcommand{\n}{\bf n}
\newcommand{\nno}{\nonumber}
\newcommand{\nord}{\mbox{\scriptsize ${\circ\atop\circ}$}}
\newtheorem{thm}{Theorem}[section]
\newtheorem{prop}[thm]{Proposition}
\newtheorem{coro}[thm]{Corollary}
\newtheorem{conj}[thm]{Conjecture}
\newtheorem{example}[thm]{Example}
\newtheorem{lem}[thm]{Lemma}
\newtheorem{rem}[thm]{Remark}
\newtheorem{de}[thm]{Definition}
\newtheorem{hy}[thm]{Hypothesis}
\makeatletter \@addtoreset{equation}{section}
\def\theequation{\thesection.\arabic{equation}}
\makeatother \makeatletter

\begin{center}
{\Large \bf On vertex Leibniz algebras}
\end{center}

\begin{center}
{Haisheng Li$^{a}$\footnote{Partially supported by NSA grant
H98230-11-1-0161 and China NSF grant (No. 11128103)}, Shaobin Tan$^{b}$\footnote{Partially supported by
China NSF grant (No.10931006) and a grant from the PhD Programs
Foundation of Ministry of Education of China (No.20100121110014).},
and Qing Wang$^{b}$\footnote{Partially supported by China NSF grant
(No.11001229) and the Fundamental Research Funds for the Central
  University (No.2012121004).}\\
$\mbox{}^{a}$Department of Mathematical Sciences\\
Rutgers University, Camden, NJ 08102, USA\\
$\mbox{}^{b}$School of Mathematical Sciences\\ Xiamen University,
Xiamen 361005, China}
\end{center}

\begin{abstract}
In this paper, we study a notion of what we call vertex Leibniz algebra.
This notion naturally extends that of vertex algebra without vacuum,
which was previously introduced by Huang and Lepowsky.
We show that every vertex algebra without vacuum
can be naturally extended to a vertex algebra. On the other hand,
we show that a vertex Leibniz algebra
can be embedded into a vertex algebra if and only if it admits a faithful module.
To each vertex Leibniz algebra we associate a vertex algebra without vacuum
which is universal to the forgetful functor.
Furthermore,  from any Leibniz algebra $\g$ we construct a vertex Leibniz algebra $V_{\g}$ and
show that $V_{\g}$ can be embedded into a vertex algebra if and only if $\g$ is a Lie algebra.
\end{abstract}

\section{Introduction}
The notion of vertex algebra has multi-faces. As one of many equivalent definitions, the notion of vertex
algebra can be defined by using existence of a vacuum and a stringy
Jacobi identity as the main axioms, where a vacuum for a vertex
algebra is analogous to an identity for a ring. This very Jacobi
identity is also known to be equivalent to what were called weak
commutativity (namely, locality) and weak associativity. This demonstrates a striking
similarity between the notions of vertex algebra and
commutative associative algebra with unit. On the other
hand, a fundamental property for vertex algebras is what is called
skew symmetry (see \cite{fhl}), which is a stringy analog of the skew symmetry for a Lie
algebra. This shows that vertex algebras are
also analogous to Lie algebras.

In literature, besides the formal variable approach (see \cite{flm}, \cite{fhl}) to vertex
algebras, there is also a $\D$-module approach (see \cite{bd}).
To bridge the $\D$-module approach and the formal variable approach,
Huang and Lepowsky (see \cite{hl}) introduced and studied
a notion of what was called a vertex algebra without vacuum.
The main axioms for a vertex algebra without vacuum in
this sense are Jacobi identity and skew symmetry (and one property
involving the translation operator).
It is clear that the notion of vertex algebra without vacuum
is a closer analog of the notion of Lie algebra (than that of vertex algebra).

In this paper, we study a notion of what we call vertex Leibniz algebra, which
naturally generalizes the notion of vertex algebra without vacuum of Huang and Lepowsky.
As it was mentioned before, one equivalent
definition of a vertex algebra is to use the existence of a vacuum
and Jacobi identity as the main axioms. The notion of vertex Leibniz
algebra is defined by simply dropping the existence of a vacuum from
this definition, to leave Jacobi identity only as the main axiom.
This notion generalizes that of vertex algebra without vacuum
in the same way that the notion of Leibniz algebra generalizes that of Lie algebra.

One of the main motivations for studying vertex Leibniz algebras came from a
previous work \cite{ltw2}, in which a theory of what we
called toroidal vertex algebras was developed. Toroidal vertex algebras are high
dimension analogues of vertex algebras to a certain extent. For a
toroidal vertex algebra, we have a Jacobi identity as the main
axiom, but we do not have a vacuum with the usual creation property,
neither do we have a skew symmetry. Another motivation came from certain studies
(\cite{loday}, \cite{loday-pir}, \cite{kw}; cf. \cite{cha}, \cite{cardona}, \cite{uc}),
in which Leibniz algebras were found to have important applications.

In this paper, we present various axiomatic results, to lay the foundation for further studies.
We first revisit Huang and Lepowsky's notion of vertex algebra without vacuum.
As an analog of a canonical exercise in classical ring theory, we show that
every vertex algebra without vacuum can be embedded into a vertex
algebra and every vertex algebra without vacuum admits a faithful
representation. (Presumably, this was known to Huang and Lepowsky.)
This shows that vertex algebras without vacuum are not much
different from vertex algebras (with vacuum). On the other hand,
vertex Leibniz algebras are quite different. From definition,
any closed subspace of a vertex algebra
in the obvious sense is automatically a vertex Leibniz algebra, but {\em not}
every vertex Leibniz algebra can be embedded into a vertex
algebra.
It is proved that a vertex Leibniz algebra $V$ can be embedded into a vertex algebra
if and only if $V$ admits a faithful representation.
This particularly implies that if the vertex operator map $Y(\cdot,x)$ of a vertex Leibniz algebra $V$
is injective, then $V$ is isomorphic to a vertex Leibniz subalgebra of
a vertex algebra.
In this paper,  to each vertex
Leibniz algebra we construct a vertex algebra without vacuum, which is universal
to the obvious forgetful functor.
Philosophically speaking, vertex Leibniz algebras are extensions of vertex algebras,
whereas vertex algebras without vacuum are ``subalgebras'' of vertex algebras.

Also among the results of this paper,  we establish an existence theorem of a vertex Leibniz algebra structure
 on a vector space with certain partial structures,
which can be considered as an analog of a theorem of Xu (see \cite{xu1}, \cite{xu}).
 As an application we construct a vertex Leibniz
algebra from any classical Leibniz algebra. More specifically, let $\g$ be a Leibniz
algebra. The tensor product $L(\g)=\g\otimes \C[t,t^{-1}]$ (with
$\C[t,t^{-1}]$ viewed as a commutative associative algebra) is a
Leibniz algebra. Extend the left action of $\g$ on $\g$ to obtain a
$\g\otimes \C[t]$-module by letting $\g\otimes t\C[t]$ act
trivially. Denote by $V_{\g}$ the $L(\g)$-module induced from the
$\g\otimes \C[t]$-module $\g$. It is proved that there exists a
canonical vertex Leibniz algebra structure on $V_{\g}$, which
contains $\g$ as a generating subspace.
Furthermore, it is shown that if $\g$ is not a Lie algebra, $V_{\g}$ is  an authentic vertex Leibniz algebra;
it cannot be embedded into a vertex algebra.

This paper is organized as follows: In Section 2, we first recall or
introduce the basic concepts of vertex algebra, vertex algebra
without vacuum, and vertex Leibniz algebra, and then we give natural
connections among them. In Section 3, we
establish an analog of a theorem of Xu and then use this to associate vertex Leibniz
algebras and their modules to loop Leibniz algebras.

{\bf Note added:} After this paper had been completed, one of us (H.L)
received a recent preprint \cite{robinson2} of Thomas Robinson, and meanwhile
Robinson's  earlier paper \cite{robinson1} was also brought to our attention.
Among other results, Robinson had studied in \cite{robinson1} a notion of vacuum free vertex algebra,
where a vacuum free vertex algebra is a vertex Leibniz algebra in the sense of this present paper such that the vertex operator map is injective,
and showed (Remark 4.1 therein) that any vacuum free vertex algebra can be embedded into a vertex algebra.
There is a certain amount of overlap between \cite{robinson2} and the present paper.
More specifically, an ertex algebra in the sense of \cite{robinson2} is the same
as a vertex Leibniz algebra, while a $D$-ertex algebra therein
amounts to a vertex algebra without vacuum in the sense of Huang-Lepowsky. Furthermore, Lemma \ref{ltwo-D-properties}
and Theorem \ref{tvano-va}, and certain results analogous to Propositions \ref{pleibniz-vano} and \ref{pembedding} of this present paper
have also been obtained in \cite{robinson2}.
We thank James Lepowsky for the information and thank Robinson for kindly sharing his preprint \cite{robinson2} with us.

\section{Vertex algebras, vertex algebras without vacuum, and vertex Leibniz algebras}
In this section, we first recall the notion of vertex algebra and
the basic properties, and we then recall the notion of vertex
algebra without vacuum from \cite{hl} and formulate a notion of
vertex Leibniz algebra. As the main results, we establish canonical
connections among these notions. For this paper, all vector spaces
are assumed to be over $\C$ (the field of complex numbers), and $\N$
denotes the set of nonnegative integers.

We start with the notion of vertex algebra.
The following is one of many equivalent definitions (\cite{bor}, \cite{flm}; cf. \cite{ll}):

\bd{dva} {\em A {\em vertex algebra} is a vector space $V$ equipped
with a linear map
\begin{eqnarray}\label{evo-map}
 Y(\cdot,x):&& V\rightarrow \Hom
(V,V((x)))\subset (\End V)[[x,x^{-1}]],\nonumber\\
&&v\mapsto Y(v,x)=\sum_{n\in \Z}v_{n}x^{-n-1}\ \ (\mbox{where
}v_{n}\in \End V)
\end{eqnarray}
and equipped with a distinguished vector ${\bf 1}\in V$, called the
{\em vacuum vector}, such that all the following conditions are
satisfied for $u,v\in V$:
$$Y({\bf 1},x)v=v\ \ (\mbox{the {\em vacuum property}}),$$
$$Y(v,x){\bf 1}\in V[[x]]\ \ \mbox{ and }\ \ \lim_{x\rightarrow
0}Y(v,x){\bf 1}=v\ \ (\mbox{the {\em creation property}}),$$
\begin{eqnarray}\label{ejacobi-def}
&&x_{0}^{-1}\delta\left(\frac{x_{1}-x_{2}}{x_{0}}\right)Y(u,x_{1})Y(v,x_{2})
-x_{0}^{-1}\delta\left(\frac{x_{2}-x_{1}}{-x_{0}}\right)Y(v,x_{2})Y(u,x_{1})\nonumber\\
&&\hspace{2cm}=x_{2}^{-1}\delta\left(\frac{x_{1}-x_{0}}{x_{2}}\right)Y(Y(u,x_{0})v,x_{2})
\end{eqnarray}
(the {\em Jacobi identity}). The linear map $Y(\cdot,x)$ is called the {\em
vertex operator map}, while for $v\in V$, $Y(v,x)$ is called the
{\em vertex operator associated to $v$}.} \ed

Let $V$ be a vertex algebra. A (resp. {\em faithful}) {\em $V$-module} is a
vector space $W$ equipped with a (resp. injective) linear map $$Y_{W}(\cdot,x):
V\rightarrow \Hom (W,W((x)))\subset (\End W)[[x,x^{-1}]]$$ such that
$Y_{W}({\bf 1},x)=1_{W}$ (the identity operator on $W$) and such that Jacobi identity holds.

Note that for a vertex algebra $V$,  due to the creation property the vertex
operator map $Y(\cdot,x)$ is {\em always} injective.
That is, the adjoint module (or representation) is always faithful.

\br{rleft-right-identity}
{\em It was proved in \cite{ll} that the vacuum property
(which is an analog of the left identity property for rings) in the definition of
the notion of vertex algebra follows from the other axioms.
This is not the case for the creation property
(which is an analog of the right identity property for rings), though
the creation property can be replaced with the condition that
the vertex operator map is injective.}
\er

Let $V$ be a vertex algebra. Define a linear operator $\D$ on $V$ by
\begin{eqnarray}\label{edef-D}
\D(v)=v_{-2}{\bf 1}=\frac{d}{dx}\left(Y(v,x){\bf 1}\right)|_{x=0}
\end{eqnarray}
for $v\in V$. Then
\begin{eqnarray}
&&[\D,Y(v,x)]=\frac{d}{dx}Y(v,x)\ \ \ \ (\mbox{the {\em
$\D$-bracket-derivative property}}),
\ \ \ \ \label{dbracket}\\
&&Y(\D v,x) =\frac{d}{dx}Y(v,x)\ \ \mbox{ for }v\in
V\label{eD-derivative}
\end{eqnarray}
(the  {\em $\D$-derivative property}), and
\begin{eqnarray}\label{eskew-symmetry}
Y(u,x)v=e^{x\D}Y(v,-x)u \ \ \mbox{ for }u,v\in V
\end{eqnarray}
(the {\em skew symmetry}).

The following is a very useful alternative definition of a vertex
algebra (see \cite{fhl}, \cite{dl}, \cite{li-local}; cf. \cite{ll}):

\bt{tdef-locality-D} A vertex algebra can be equivalently defined as
a vector space $V$ equipped with a linear map $Y(\cdot,x)$ as in
(\ref{evo-map}), a vector ${\bf 1}\in V$, satisfying the vacuum property and creation property,
and equipped with a linear operator $D$ on $V$, satisfying $D{\bf 1}=0$, the
$D$-bracket-derivative property and the {\em locality}: For any
$u,v\in V$, there exists a nonnegative integer $k$ such that
\begin{eqnarray}
(x_{1}-x_{2})^{k}[Y(u,x_{1}),Y(v,x_{2})]=0.
\end{eqnarray}
\et

The following notion is due to Huang and Lepowsky (see \cite{hl}):

\bd{dvawv} {\em A {\em vertex algebra without vacuum} is a vector
space $V$ equipped with a linear map $Y(\cdot,x)$ as in (\ref{evo-map}) and a
linear operator $D$ on $V$ such that Jacobi identity, the
$D$-derivative property (\ref{eD-derivative}), and skew symmetry
(\ref{eskew-symmetry}) hold. } \ed

\br{rvawv-equivdef}
{\em It was known (see \cite{fhl}) that  in the presence of skew symmetry,
locality, namely weak commutativity, is equivalent to weak associativity.
Therefore, in the definition of a vertex algebra without vacuum, Jacobi identity can be equivalently
replaced by either locality or weak associativity.}
\er

\br{rborcherds} {\em Let $A$ be a commutative associative algebra (possibly without a unit)
with a derivation $D$.
For $a\in A$, define $Y(a,x)\in (\End_{\C} A)[[x]]$ by $Y(a,x)b=(e^{xD}a)b$ for $b\in A$.
Then $(A,Y,D)$ carries the structure of a vertex algebra without vacuum (see \cite{bor}).
 On the other hand, assume that $(V,Y,D)$ is a vertex algebra without vacuum such that
$Y(u,x)v\in V[[x]]$ for $u,v\in V$. Then one can show that $V$ is a commutative associative algebra
with $D$ as a derivation, where $u\cdot v=Y(u,x)v|_{x=0}$ for $u,v\in V$.
Furthermore, we have $Y(u,x)v=(e^{xD}u)\cdot v$ for $u,v\in V$.}
\er

\bl{ltwo-D-properties} Let $V$ be
a vector space equipped with a linear map $Y(\cdot,x)$ as in
(\ref{evo-map}) and let $D$ be a linear operator on $V$ such that
skew symmetry holds. Then the $D$-bracket-derivative property is
equivalent to the $D$-derivative property. \el

\begin{proof} For $u,v\in V$, using skew symmetry one gets
\begin{eqnarray*}
&&Y(Du,x)v-\frac{d}{dx}Y(u,x)v\\
&=&e^{xD}Y(v,-x)Du-e^{xD}DY(v,-x)u-e^{xD}\frac{d}{dx}Y(v,-x)u\\
&=&e^{xD}\left(-[D,Y(v,-x)]+\frac{d}{d(-x)}Y(v,-x)\right)u.
\end{eqnarray*}
Then the equivalence follows immediately.
\end{proof}

Define a notion of homomorphism of vertex algebras without vacuum
in the obvious way: A {\em homomorphism} of vertex
algebras without vacuum from $U$ to $V$ is a linear map $\psi:
U\rightarrow V$, satisfying the condition that
\begin{eqnarray}
\psi (Du)=D\psi(u) \ \ \mbox{ and }\ \
\psi(Y(u,x)u')=Y(\psi(u),x)\psi(u')
\end{eqnarray}
for $u,u'\in U$.

We have the following analog of a well known classical result in ring
theory:

\bt{tvano-va} Let $(V,Y,D)$ be a vertex algebra without vacuum. Set
$\overline{V}=V\oplus \C{\bf 1}$, where $\C {\bf 1}$ is a
one-dimensional vector space with a distinguished base vector ${\bf
1}$. Define a linear map
$$\overline{Y}(\cdot,x):
\overline{V}\rightarrow (\End \overline{V})[[x,x^{-1}]]$$ by
$\overline{Y}({\bf 1},x)=1$ and
$$\overline{Y}(u,x)(v+\lambda{\bf 1})=Y(u,x)v+\lambda e^{xD}u$$
for $u,v\in V,\ \lambda\in \C$. Then
$(\overline{V},\overline{Y},{\bf 1})$ carries the structure of a
vertex algebra with $V$ as a vertex subalgebra without vacuum. \et

\begin{proof} It is clear that
$\overline{Y}(\cdot,x)$ maps $\overline{V}$ to $\Hom
(\overline{V},\overline{V}((x)))$. Extend $D$ to a linear operator
$\overline{D}$ on $\overline{V}$ by defining $\overline{D}({\bf
1})=0$.
 For $u\in V,\ \alpha\in \C$, by definition we have
$$\overline{Y}(u+\alpha {\bf 1},x){\bf
1}=e^{xD}u+\alpha {\bf 1}=e^{x\overline{D}}(u+\alpha {\bf 1}).$$ We
see that the vacuum property and the creation property hold.

For $u,v\in V,\ \alpha\in \C$, using $D$-bracket-derivative property
(recall Lemma \ref{ltwo-D-properties}) we get
\begin{eqnarray*}
\overline{D}\cdot \overline{Y}(u,x)(v+\alpha{\bf 1})
&=&D(Y(u,x)v+\alpha e^{xD}u)\\
&=&Y(u,x)Dv+\frac{d}{dx}Y(u,x)v+\alpha De^{xD}u\\
&=&\overline{Y}(u,x)\overline{D}(v+\alpha{\bf
1})+\frac{d}{dx}\overline{Y}(u,x)(v+\alpha {\bf 1}).
\end{eqnarray*}
This proves
$$[\overline{D},\overline{Y}(u,x)]=\frac{d}{dx}\overline{Y}(u,x)\ \ \mbox{ for }u\in V.$$
We also have $[\overline{D},\overline{Y}({\bf
1},x)]=\frac{d}{dx}\overline{Y}({\bf 1},x)$.

Let $u,v,w\in V,\; \alpha \in \C$. By definition we have
\begin{eqnarray*}
\overline{Y}(u,x_{1})\overline{Y}(v,x_{2})(w+\alpha{\bf 1})
&=&\overline{Y}(u,x_{1})(Y(v,x_{2})w+\alpha e^{x_{2}D}v)\\
&=&Y(u,x_{1})Y(v,x_{2})w+\alpha Y(u,x_{1})e^{x_{2}D}v,
\end{eqnarray*}
\begin{eqnarray*}
\overline{Y}(v,x_{2})\overline{Y}(u,x_{1})(w+\alpha {\bf
1})=Y(v,x_{2})Y(u,x_{1})w+\alpha Y(v,x_{2})e^{x_{1}D}u.
\end{eqnarray*}
Furthermore, using $D$-bracket-derivative property and skew symmetry
we get
\begin{eqnarray*}
&&Y(u,x_{1})e^{x_{2}D}v=e^{x_{2}D}Y(u,x_{1}-x_{2})v=e^{x_{1}D}Y(v,-x_{1}+x_{2})u,\\
&&Y(v,x_{2})e^{x_{1}D}u=e^{x_{1}D}Y(v,x_{2}-x_{1})u.
\end{eqnarray*}
Let $k$ be a nonnegative integer such that
$$(x_{1}-x_{2})^{k}[Y(u,x_{1}),Y(v,x_{2})]=0
\ \ \mbox{ and } \ \ x^{k}Y(v,x)u\in V[[x]],$$ so that
$$(x_{1}-x_{2})^{k}Y(v,x_{2}-x_{1})u=(x_{1}-x_{2})^{k}Y(v,-x_{1}+x_{2})u.$$
Then we obtain
\begin{eqnarray*}
(x_{1}-x_{2})^{k}\overline{Y}(u,x_{1})\overline{Y}(v,x_{2})
=(x_{1}-x_{2})^{k}\overline{Y}(v,x_{2})\overline{Y}(u,x_{1}).
\end{eqnarray*}
By Theorem \ref{tdef-locality-D}, $(\overline{V},\overline{Y},{\bf
1})$ carries the structure of a vertex algebra. It can be readily
seen that $V$ is a vertex subalgebra without vacuum of
$\overline{V}$.
\end{proof}

\bd{dmodule-vawv} {\em Let $(V,Y,D)$ be a vertex algebra without vacuum. A
(resp. {\em faithful}) {\em $V$-module} is a vector space $W$
equipped with a (resp. injective) linear map
$$Y_{W}(\cdot,x): V\rightarrow \Hom (W,W((x)))\subset (\End
W)[[x,x^{-1}]]$$ such that Jacobi identity holds and such that
\begin{eqnarray}\label{ed-module}
Y_{W}(D v,x)=\frac{d}{dx}Y_{W}(v,x) \ \ \ \mbox{ for }v\in V.
\end{eqnarray}}
\ed

As an immediate consequence of Theorem \ref{tvano-va}, we have:

\bc{cfaithful-module} Every vertex algebra without vacuum admits a
faithful module. \ec

The following shows that for a vertex algebra without vacuum $V$, a $V$-module
exactly amounts to a $\overline{V}$-module:

\bp{pmodule-vawv}
Let $V$ be a vertex algebra without vacuum
and let $\overline{V}=V\oplus \C {\bf 1}$ be the vertex algebra
constructed in Theorem \ref{tvano-va}. Let $(W,Y_{W})$ be any $V$-module.
Extend $Y_{W}(\cdot,x)$ to a linear map
$$\overline{Y}_{W}(\cdot,x):\ \overline{V}\rightarrow (\End W)[[x,x^{-1}]]$$
by defining $\overline{Y}_{W}({\bf 1},x)=1_{W}$. Then
$(W,\overline{Y}_{W})$ is a $\overline{V}$-module.
 \ep

\begin{proof} We need to prove that the Jacobi identity for any pair $(u,v)$ in $V\cup \{\bf 1\}$ holds. If
$u,v\in V$, or $u=v={\bf 1}$, it is clear. If $u={\bf 1}, v\in V$,
it is also clear. Consider the case with $u\in V,\ v={\bf 1}$. Using
(\ref{ed-module}) we have
$$Y_{W}(Y(u,x_{0}){\bf
1},x_{2})=Y_{W}(e^{x_{0}D}u,x_{2})=Y_{W}(u,x_{2}+x_{0}).$$ Then the
Jacobi identity for the pair $(u,{\bf 1})$ follows from basic
delta-function properties.
\end{proof}

Next, we study an analog of the notion of Leibniz algebra.

\bd{dleibniz} {\em A {\em vertex Leibniz algebra} is a vector space
$V$ equipped with a linear map $Y(\cdot,x)$ as in (\ref{evo-map})
such that Jacobi identity (\ref{ejacobi-def}) holds.} \ed

\br{rborcherds-leibniz} {\em
Let $A$ be an associative algebra (possibly without a unit)
such that $a(bc)=b(ac)$ for $a,b,c\in A$.
Such an algebra is called a {\em Perm algebra} in literature (see \cite{cha}, \cite{cl}).
Let $D$ be a linear endomorphism of $A$ such that
$$\left(D(ab)-D(a)b-a D(b)\right)c=0\ \ \ \mbox{ for }a,b,c\in A.$$
For $a\in A$, define $Y(a,x)\in (\End_{\C} A)[[x]]$ by $Y(a,x)b=(e^{xD}a)b$ for $b\in A$.
It is straightforward to show that $(A,Y)$ carries the structure of a vertex Leibniz algebra such that
 the $D$-derivative property (\ref{eD-derivative}) holds.
On the other hand, assume that $(V,Y)$ is a vertex Leibniz algebra  such that
$Y(u,x)v\in V[[x]]$ for $u,v\in V$. Define $u\cdot v=(Y(u,x)v)|_{x=0}$ for $u,v\in V$.
It is straightforward to show that $(V,\cdot)$ is a Perm algebra. If we in addition assume that
there exists a linear operator $D$ on $V$ such that
 the $D$-derivative property (\ref{eD-derivative}) holds, then
 $Y(u,x)v=(e^{xD}u)\cdot v$ for $u,v\in V$. Furthermore, for $u,v,w\in V$,
 as $Y(u,x_{1})w\in V[[x_{1}]]$, applying $\Res_{x_{1}}$ to the Jacobi identity
 for the triple $(u,v,w)$ we get
 $$Y(Y(u,x_{0})v,x_{2})w=Y(u,x_{0}+x_{2})Y(v,x_{2})w.$$
 Extracting the coefficient of $x_{0}^{0}x_{2}$ from both sides we obtain
 $$\left(D(u\cdot v)-D(u)\cdot v-u\cdot D(v)\right)\cdot w=0\ \ \ \mbox{ for }u,v,w\in V.$$}
\er

\br{rclosed-va} {\em By definition, a vertex algebra (equipped with the canonical
operator $\D$) is a vertex algebra without vacuum, and a vertex
algebra without vacuum is a vertex Leibniz algebra.
On the other hand, let $V$ be a vertex algebra and suppose that $K$
is a subspace which is {\em closed } in the sense that
$$a_{n}b\in K\ \ \mbox{ for }a,b\in K,\ n\in \Z.$$
It is clear that $K$ is a vertex Leibniz algebra. Furthermore, if
$K$ is also $\D$-stable (recall (\ref{edef-D})), then $K$ equipped with $\D$ (viewed as an
operator on $K$) is a vertex algebra without vacuum.} \er

The following is a straightforward analog of a classical fact in the theory of
Leibniz algebras (cf. \cite{cardona}):

\bl{lhom-construction}
Let $V$ be a vertex Leibniz algebra and let $W$
be a $V$-module with a $V$-module homomorphism $\psi: W\rightarrow V$.
Define a linear map
$Y(\cdot,x) : W\rightarrow (\End W)[[x,x^{-1}]]$ by
$$Y(w,x)w'=Y(\psi(w),x)w'\ \ \ \mbox{ for }w,w'\in W.$$
Then $(W,Y)$ carries the structure of a vertex Leibniz algebra and
$\psi$ is a homomorphism of vertex Leibniz algebras.
\el

\bex{eaxmple-verma}
{\em Let $\g$ be a general Lie algebra equipped with a
non-degenerate symmetric invariant bilinear form $\<\cdot,\cdot\>$.
Associated to $(\g, \<\cdot,\cdot\>)$, one has an affine Lie algebra
$\hat{\g}=\g\otimes \C[t,t^{-1}]\oplus \C {\bf k}$.
For any complex number $\ell$, we have a vertex algebra
$V_{\hat{\g}}(\ell,0)$ (cf. \cite{fz}), where
$$V_{\hat{\g}}(\ell,0)=U(\hat{\g})\otimes _{U(\g\otimes \C[t]+\C {\bf k})}\C_{\ell},$$
a generalized Verma $\hat{\g}$-module induced from a $1$-dimensional $(\g\otimes \C[t]+\C {\bf k})$-module $\C_{\ell}=\C$
on which $\g\otimes \C[t]$ acts trivially and ${\bf k}$ acts as scalar $\ell$.
It is also known (cf. \cite{li-local}) that every restricted $\hat{\g}$-module $W$ of level $\ell$
is naturally a $V_{\hat{\g}}(\ell,0)$-module.
Assume that $\g$ is a finite-dimensional simple Lie algebra. Denote by $M_{\hat{\g}}(\ell,0)$
the Verma $\hat{\g}$-module of level $\ell$ with highest weight $0$. Then
$M_{\hat{\g}}(\ell,0)$ is a $V_{\hat{\g}}(\ell,0)$-module and the natural $\hat{\g}$-module homomorphism from
$M_{\hat{\g}}(\ell,0)$ to $V_{\hat{\g}}(\ell,0)$  is a $V_{\hat{\g}}(\ell,0)$-module homomorphism.
By Lemma \ref{lhom-construction}, $M_{\hat{\g}}(\ell,0)$ is naturally
a vertex Leibniz algebra.}
\eex

\br{rmodel-construction} {\em Note that in Lemma \ref{lhom-construction},
for $w\in\ker \psi, \ w'\in W$ we have $Y(w,x)w'=0$.
Thus, if $\psi$ is not one-to-one, $W$ as a vertex Leibniz algebra
cannot have a vacuum vector so that $W$ cannot be a vertex algebra
even if $V$ is a vertex algebra.}
\er

\bex{kwconstruction} {\em A special case of Lemma \ref{lhom-construction} gives an analog of
 the {\em hemisemidirect product} $\g \ltimes_{H}W$ of a Lie algebra $\g$ with a $\g$-module $W$
(cf. \cite{kw}).
Let $V$ be a vertex algebra and let $W$ be a $V$-module. Then $W\oplus V$ is a $V$-module where the projection
of $V\oplus W$ onto $V$ is a $V$-module homomorphism. In view of Lemma \ref{lhom-construction},
$V+W$ becomes a vertex Leibniz algebra with
\begin{eqnarray}
Y(u+w,x)(v+w')=Y(u,x)v+Y(u,x)w'\ \ \ \mbox{ for }u,v\in V,\ w,w'\in W.
\end{eqnarray}
{}From Remark \ref{rmodel-construction}, if $W\ne 0$, $V+W$ is not a vertex algebra.
 In fact, if $W\ne 0$, $V+W$ even cannot be embedded
into a vertex algebra. This is because if $V+W$ is a vertex Leibniz subalgebra of a vertex algebra $K$, then
$$w=Y({\bf 1},x)w= e^{x\D}Y(w,-x){\bf 1}=0\ \ \mbox{ for all }w\in W,$$
where ${\bf 1}$ is the vacuum vector of $V$ and $\D$ is the $\D$-operator of $K$.} \eex

\bex{example-orbifold} {\em Let $V$ be a vertex algebra and let $G$
be a finite automorphism group of $V$.  Denote by $V^{G}$ the vertex subalgebra of $G$-fixed points in $V$
and define a linear map $\psi: V\rightarrow V^{G}$ by
$$\psi(v)=\frac{1}{|G|}\sum_{g\in G}gv\ \ \ \mbox{ for }v\in V.$$
It can be readily seen that $\psi$ is a $V^{G}$-module homomorphism.
In view of Lemma \ref{lhom-construction}, $V$ becomes a vertex Leibniz algebra with the vertex operator map
$Y_{G}(\cdot,x): V\rightarrow (\End V)[[x,x^{-1}]]$ given by
\begin{eqnarray}
Y_{G}(u,x)v=\frac{1}{|G|}\sum_{g\in G}Y(gu,x)v
\end{eqnarray}
for $u,v\in V$. We have
\begin{eqnarray*}
Y_{G}(hu,x)v=Y_{G}(u,x)v,\ \ \ \ Y_{G}(u,x)(hv)=hY_{G}(u,x)v
\end{eqnarray*}
for $h\in G,\ u,v\in V$. On the other hand, by Lemma \ref{lhom-construction},
 $\psi$ is a vertex Leibniz algebra homomorphism.
It is straightforward to show that $\ker \psi$ is exactly
the linear span of $gu-u$ for $g\in G,\ u\in V$.
Consequently, $\psi$ gives rise to an isomorphism from the quotient vertex Leibniz algebra
$V/\ker \psi$ to vertex algebra $V^{G}$.  } \eex

\bp{pvla-vawv} Let $V$ be a vertex Leibniz algebra and let $D$ be a
linear operator on $V$ such that the $D$-derivative property holds.
If $Y(\cdot,x)$ is injective, then $V$ is a vertex algebra without
vacuum. \ep

\begin{proof} We need to show that skew symmetry holds. This basically follows from \cite{fhl}.
Let $u,v,w\in V$.
Note that the left-hand side of Jacobi identity has a symmetry
$$(u,v,x_{1},x_{2},x_{0})\leftrightarrow (v,u,x_{2},x_{1},-x_{0}),$$
 from which one obtains
\begin{eqnarray}\label{epre-skew-symmetry}
Y(Y(u,x_{0})v,x_{2})w=Y(Y(v,-x_{0})u,x_{2}+x_{0})w.
\end{eqnarray}
Using this and the $D$-derivative property we get
\begin{eqnarray}
Y(Y(u,x_{0})v,x_{2})=Y(Y(v,-x_{0})u,x_{2}+x_{0})=Y(e^{x_{0}D}Y(v,-x_{0})u,x_{2}).
\end{eqnarray}
Since $Y(\cdot,x)$ is injective, we immediately obtain skew
symmetry.
\end{proof}

Note that for a vertex algebra $V$, the operator $D$ is canonically
defined by $D(v)=v_{-2}{\bf 1}$ for $v\in V$. Next, we discuss the
uniqueness of operator $D$ for a vertex algebra without vacuum.

Let $V$ be a vertex Leibniz algebra. Set
\begin{eqnarray}
I_{V}=\{ v\in V\ |\ Y(v,x)=0\}.
\end{eqnarray}
Then $I_{V}$ is an {\em ideal} of $V$ in the sense that $I_{V}$ is a
subspace satisfying
$$u_{n}v,\ \ v_{n}u\in I_{V}\ \ \ \mbox{ for all }u\in V,\ v\in I_{V}.$$
Let $\Der (V,I_{V})$ denote the set of linear maps $\psi\in
\Hom(V,I_{V})$, satisfying
\begin{eqnarray}
\psi(Y(u,x)v)=Y(\psi(u),x)v+Y(u,x)\psi(v)\ \ \ \mbox{ for }u,v\in V.
\end{eqnarray}
We have:

\bp{puniqueness-D} Let $V$ be a vertex Leibniz algebra. Suppose that
$D_{0}$ is a linear operator on $V$, satisfying (\ref{dbracket}) and
(\ref{eD-derivative}). Then $D_{0}+\Der (V,I_{V})$ is exactly the
set of linear operators $D$ satisfying (\ref{dbracket}) and (\ref{eD-derivative}). In
particular, if $V$ is a vertex algebra without vacuum such that $Y(\cdot,x)$ is injective, then the
operator $D$ is unique.\ep

\begin{proof} Let $\psi\in \Der (V,I_{V})$.
For $u,v\in V$, since $Y(\psi(u),x)=0$ we have
$$\psi(Y(u,x)v)=Y(\psi(u),x)v+Y(u,x)\psi(v)=Y(u,x)\psi(v).$$
Then
\begin{eqnarray*}
&&[D_{0}+\psi,Y(u,x)]=[D_{0},Y(u,x)]=\frac{d}{dx}Y(u,x),\\
&&Y((D_{0}+\psi)u,x)=Y(D_{0}u,x)=\frac{d}{dx}Y(u,x).
\end{eqnarray*}
On the other hand, suppose that $D$ is a linear operator on $V$,
satisfying (\ref{dbracket}) and (\ref{eD-derivative}). Set
$\psi=D-D_{0}$. For $u\in V$, we have
$$Y(\psi(u),x)=Y(Du,x)-Y(D_{0}u,x)=\frac{d}{dx}Y(u,x)-\frac{d}{dx}Y(u,x)=0,$$
which implies $\psi(u)\in I_{V}$. Furthermore,
$$[\psi,Y(u,x)]=[D,Y(u,x)]-[D_{0},Y(u,x)]=\frac{d}{dx}Y(u,x)-\frac{d}{dx}Y(u,x)=0.$$
Thus $\psi\in \Der (V,I_{V})$.
\end{proof}

\bd{dmodule-lva} {\em Let $V$ be a vertex Leibniz algebra. A (resp.
{\em faithful}) {\em $V$-module} is a vector space $W$ equipped with
a (resp. injective) linear map
$$Y_{W}(\cdot,x):
V\rightarrow \Hom (W,W((x)))\subset (\End W)[[x,x^{-1}]]$$ such that
Jacobi identity holds.} \ed

We have (see also \cite{robinson1}, Remark 4.1):

\bp{pfaithful} Let $V$ be a vertex Leibniz algebra. Then $V$ admits a
faithful module if and only if $V$ is a vertex Leibniz subalgebra of some vertex algebra.
In particular, if the vertex
operator map $Y(\cdot,x)$ of $V$ is injective, i.e.,
for any nonzero $u\in V$, there exists $v\in V$ such that $Y(u,x)v\ne 0$, then $V$ is a vertex Leibniz
subalgebra of some vertex algebra.
\ep

\begin{proof} Since the adjoint module for every vertex algebra is faithful, the
``if'' part is clear. Now, assume that $V$ has a faithful module
$(W,Y_{W})$. Set $U=\{ Y_{W}(v,x)\ |\ v\in V\}$. It follows from
Jacobi identity that $U$ is a local subset of $\Hom
(W,W((x)))$. From \cite{li-local}, $U$ generates a vertex algebra
$\<U\>$ inside $\Hom (W,W((x)))$, where the identity operator $1_{W}$ on $W$ serves as the vacuum vector
and the vertex operator map, denoted by $Y_{\E}(\cdot,x)$, is given by
\begin{eqnarray*}
&&Y_{\E}(\alpha(x),x_{0})\beta(x)\\
&=&\Res_{x_{1}}\left(x_{0}^{-1}\delta\left(\frac{x_{1}-x}{x_{0}}\right)\alpha(x_{1})\beta(x)
-x_{0}^{-1}\delta\left(\frac{x-x_{1}}{-x_{0}}\right)\beta(x)\alpha(x_{1})\right)
\end{eqnarray*}
for $\alpha(x),\beta(x)\in \<U\>$. For $u,v\in V$, taking $\alpha(x)=Y_{W}(u,x),\ \beta(x)=Y_{W}(v,x)$ we get
\begin{eqnarray*}
&&Y_{\E}(Y_{W}(u,x),x_{0})Y_{W}(v,x)\\
&=&\Res_{x_{1}}\left(x_{0}^{-1}\delta\left(\frac{x_{1}-x}{x_{0}}\right)Y_{W}(u,x_{1})Y_{W}(v,x)
-x_{0}^{-1}\delta\left(\frac{x-x_{1}}{-x_{0}}\right)Y_{W}(v,x)Y_{W}(u,x_{1})\right)\\
&=&Y_{W}(Y(u,x_{0})v,x).
\end{eqnarray*}
This shows that the map $Y_{W}(\cdot,x)$ from $V$ to $\<U\>$
is a homomorphism of vertex Leibniz algebras. As $Y_{W}(\cdot,x)$ is assumed to be injective,
$V$ is isomorphic to a vertex Leibniz subalgebra of $\<U\>$.
\end{proof}

Let $V$ be a vertex Leibniz algebra. Recall that $I_{V}=\{ v\in V\
|\ Y(v,x)=0\}$, which is an ideal of $V$.  One sees that $V$ is a
faithful module for the quotient vertex Leibniz algebra $V/I_{V}$.
In view of Proposition \ref{pfaithful}, $V/I_{V}$ is isomorphic to a
vertex Leibniz subalgebra of some vertex algebra. The following
generalization is also immediate:

\bc{cmodule-quotient} Let $V$ be a vertex Leibniz algebra and let
$(W,Y_{W})$ be any $V$-module. Set
$$I_{V}(W)=\{ v\in V\ |\ Y_{W}(v,x)=0\}.$$
Then $I_{V}(W)$ is an ideal of $V$ and the quotient algebra
$V/I_{V}(W)$ is isomorphic to a vertex Leibniz subalgebra of some
vertex algebra. \ec

\bd{djacobson} {\em Let $V$ be a vertex Leibniz algebra. Define
\begin{eqnarray}
J_{V}=\cap_{W\in {\mathcal{M}}(V)} \ker Y_{W}(\cdot,x),
\end{eqnarray}
an ideal of $V$, where ${\mathcal{M}}(V)$ denotes the collection of
all $V$-modules.} \ed

\bp{pJV-I} Let $V$ be a vertex Leibniz algebra.  If $I$ is an ideal
of $V$ such that $V/I$ is isomorphic to a vertex Leibniz subalgebra
of some vertex algebra, then $I\supset J_{V}$.  On the other hand,
the quotient vertex Leibniz algebra $V/J_{V}$ is isomorphic to a
vertex Leibniz subalgebra of some vertex algebra. \ep

\begin{proof} Assume that $V/I$ is a vertex Leibniz subalgebra of a vertex
algebra $K$. Then $K$ is a faithful $V/I$-module, due to the
creation property of $K$. This $V/I$-module structure on $K$ gives
rise to a $V$-module structure $Y_{K}(\cdot,x)$ with $\ker Y_{K}=I$.
Thus $J_{V}\subset \ker Y_{K}=I$. For the second assertion, note
that the direct sum of all $V$-modules is a faithful module for
$V/J_{V}$. By Proposition \ref{pfaithful}, $V/J_{V}$ is isomorphic
to a vertex Leibniz subalgebra of some vertex algebra.
\end{proof}

As an immediate consequence we have:

\bc{cJV=0}
Let $V$ be a vertex Leibniz algebra. Then $V$ can be embedded  into a
vertex algebra if and only if $J_{V}=0$.
\ec

For convenience we formulate the following notion:

\bd{dvlawithD} {\em
A {\em vertex Leibniz algebra with $D$} is a vertex Leibniz algebra $V$ equipped with a
linear operator $D$ on $V$ such that both the $D$-bracket-derivative property and the $D$-derivative
property hold.}
\ed

We have:

\bp{pleibniz-vano} Let $V$ be a vertex Leibniz algebra. Set
$$\overline{V}=\C[D]\otimes V,$$
where $D$ is a symbol.  View $D$ as a linear operator on
$\overline{V}$ in the obvious way and define a linear map
$$\overline{Y}(\cdot,x):\ \overline{V}\rightarrow (\End
\overline{V})[[x,x^{-1}]]$$ by
\begin{eqnarray}\label{edef-generating}
\overline{Y}(e^{zD}u,x)(e^{z_{0}D}v)=e^{z_{0}D}Y(u,x+z-z_{0})v
\end{eqnarray}
for $u,v\in V$. Then $(\overline{V},\overline{Y}, D)$ carries the
structure of a vertex Leibniz algebra with $D$.
Furthermore, $\overline{V}$ contains $V$ as a vertex Leibniz subalgebra.
 \ep

\begin{proof} For $u,v\in V$, applying $\frac{\partial}{\partial z_{0}}$
to (\ref{edef-generating}) we get
\begin{eqnarray*}
&&\overline{Y}(e^{zD}u,x)D(e^{z_{0}D}v)\\
&=&De^{z_{0}D}Y(u,x+z-z_{0})v+e^{z_{0}D}\frac{\partial}{\partial
z_{0}}Y(u,x+z-z_{0})v\\
&=&De^{z_{0}D}Y(u,x+z-z_{0})v-\frac{\partial}{\partial
x}e^{z_{0}D}Y(u,x+z-z_{0})v\\
&=&D\overline{Y}(e^{zD}u,x)(e^{z_{0}D}v)-\frac{\partial}{\partial
x}\overline{Y}(e^{zD}u,x)(e^{z_{0}D}v).
\end{eqnarray*}
This gives
\begin{eqnarray*}
\left[D,\overline{Y}(e^{zD}u,x)\right] =\frac{\partial}{\partial
x}\overline{Y}(e^{zD}u,x).
\end{eqnarray*}
On the other hand, applying $\frac{\partial}{\partial z}$ to
(\ref{edef-generating}) we get
\begin{eqnarray*}
\overline{Y}(De^{zD}u,x)(e^{z_{0}D}v)&=&\frac{\partial}{\partial
z}e^{z_{0}D}Y(u,x+z-z_{0})v\\
&=&\frac{\partial}{\partial x}e^{z_{0}D}Y(u,x+z-z_{0})v\\
&=&\frac{\partial}{\partial x}\overline{Y}(e^{zD}u,x)(e^{z_{0}D}v).
\end{eqnarray*}
This proves $\overline{Y}(Da,x)=\frac{d}{dx}\overline{Y}(a,x)$ for
$a\in \overline{V}$.

Let $u,v,w\in V$. By definition we have
\begin{eqnarray}
&&\overline{Y}(e^{z_{1}D}u,x_{1})\overline{Y}(e^{z_{2}D}v,x_{2})(e^{z_{0}D}w)\nonumber\\
&=&\overline{Y}(e^{z_{1}D}u,x_{1})e^{z_{0}D}Y(v,x_{2}+z_{2}-z_{0})w\nonumber\\
&=&e^{z_{0}D}Y(u,x_{1}+z_{1}-z_{0})Y(v,x_{2}+z_{2}-z_{0})w\nonumber\\
&=&e^{z_{0}D}e^{(z_{1}-z_{0})\frac{\partial}{\partial
x_{1}}}e^{(z_{2}-z_{0})\frac{\partial}{\partial
x_{2}}}Y(u,x_{1})Y(v,x_{2})w.
\end{eqnarray}
Symmetrically, we have
\begin{eqnarray}
&&\overline{Y}(e^{z_{2}D}v,x_{2})\overline{Y}(e^{z_{1}D}u,x_{1})(e^{z_{0}D}w)\nonumber\\
&=&e^{z_{0}D}e^{(z_{1}-z_{0})\frac{\partial}{\partial
x_{1}}}e^{(z_{2}-z_{0})\frac{\partial}{\partial
x_{2}}}Y(v,x_{2})Y(u,x_{1})w.
\end{eqnarray}
On the other hand, we have
\begin{eqnarray}
&&\overline{Y}\left(\overline{Y}(e^{z_{1}D}u,x_{0})(e^{z_{2}D}v),x_{2}\right)(e^{z_{0}D}w)
\nonumber\\
&=&\overline{Y}\left(e^{z_{2}D}Y(u,x_{0}+z_{1}-z_{2})v,x_{2}\right)(e^{z_{0}D}w)
\nonumber\\
&=&e^{z_{0}D}Y(Y(u,x_{0}+z_{1}-z_{2})v,x_{2}+z_{2}-z_{0})w\nonumber\\
&=&e^{z_{0}D}e^{(z_{1}-z_{2})\frac{\partial}{\partial
x_{0}}}e^{(z_{2}-z_{0})\frac{\partial}{\partial
x_{2}}}Y(Y(u,x_{0})v,x_{2})w.
\end{eqnarray}
Note that if $A(x_{0},x_{1},x_{2})$ is one of the following three
delta functions
$$x_{0}^{-1}\delta\left(\frac{x_{1}-x_{2}}{x_{0}}\right), \ \
x_{0}^{-1}\delta\left(\frac{x_{2}-x_{1}}{-x_{0}}\right),\ \
x_{1}^{-1}\delta\left(\frac{x_{2}+x_{0}}{x_{1}}\right),$$
we have
\begin{eqnarray*}
\frac{\partial}{\partial
x_{1}}A(x_{0},x_{1},x_{2})=-\frac{\partial}{\partial
x_{2}}A(x_{0},x_{1},x_{2})=-\frac{\partial}{\partial
x_{0}}A(x_{0},x_{1},x_{2}).
\end{eqnarray*}
Consequently,
\begin{eqnarray*}
e^{(z_{1}-z_{2})\frac{\partial}{\partial
x_{0}}}e^{(z_{1}-z_{0})\frac{\partial}{\partial
x_{1}}}e^{(z_{2}-z_{0})\frac{\partial}{\partial x_{2}}}
A(x_{0},x_{1},x_{2})=A(x_{0},x_{1},x_{2}).
\end{eqnarray*}
Notice that
\begin{eqnarray*}
&&e^{(z_{1}-z_{2})\frac{\partial}{\partial
x_{0}}}e^{(z_{1}-z_{0})\frac{\partial}{\partial
x_{1}}}e^{(z_{2}-z_{0})\frac{\partial}{\partial x_{2}}}
\left(Y(u,x_{1})Y(v,x_{2})w\right)\\
&=&e^{(z_{1}-z_{0})\frac{\partial}{\partial
x_{1}}}e^{(z_{2}-z_{0})\frac{\partial}{\partial x_{2}}}
\left(Y(u,x_{1})Y(v,x_{2})w\right),
\end{eqnarray*}
\begin{eqnarray*}
 &&e^{(z_{1}-z_{2})\frac{\partial}{\partial
x_{0}}}e^{(z_{1}-z_{0})\frac{\partial}{\partial
x_{1}}}e^{(z_{2}-z_{0})\frac{\partial}{\partial x_{2}}}
\left(Y(v,x_{2})Y(u,x_{1})w\right)\\
&=&e^{(z_{1}-z_{0})\frac{\partial}{\partial
x_{1}}}e^{(z_{2}-z_{0})\frac{\partial}{\partial x_{2}}}
\left(Y(v,x_{2})Y(u,x_{1})w\right),
\end{eqnarray*}
\begin{eqnarray*}
&&e^{(z_{1}-z_{2})\frac{\partial}{\partial
x_{0}}}e^{(z_{1}-z_{0})\frac{\partial}{\partial
x_{1}}}e^{(z_{2}-z_{0})\frac{\partial}{\partial
x_{2}}}\left(Y(Y(u,x_{0})v,x_{2})w\right)\\
&=&e^{(z_{1}-z_{2})\frac{\partial}{\partial
x_{0}}}e^{(z_{2}-z_{0})\frac{\partial}{\partial
x_{2}}}\left(Y(Y(u,x_{0})v,x_{2})w\right).
\end{eqnarray*}
By definition, for the given vertex Leibniz algebra $(V,Y)$, we have
the Jacobi identity for the ordered triple $(u,v,w)$. Applying
$e^{(z_{1}-z_{2})\frac{\partial}{\partial
x_{0}}}e^{(z_{1}-z_{0})\frac{\partial}{\partial
x_{1}}}e^{(z_{2}-z_{0})\frac{\partial}{\partial x_{2}}}$ to the both
sides of the Jacobi identity for the ordered triple $(u,v,w)$, we
obtain
\begin{eqnarray*}
&&x_{0}^{-1}\delta\left(\frac{x_{1}-x_{2}}{x_{0}}\right)
\overline{Y}(e^{z_{1}D}u,x_{1})\overline{Y}(e^{z_{2}D}v,x_{2})(e^{z_{0}D}w)\\
&&\ \ \ \ -x_{0}^{-1}\delta\left(\frac{x_{2}-x_{1}}{-x_{0}}\right)
\overline{Y}(e^{z_{2}D}v,x_{2})\overline{Y}(e^{z_{1}D}u,x_{1})(e^{z_{0}D}w)\\
&=&x_{2}^{-1}\delta\left(\frac{x_{1}-x_{0}}{x_{2}}\right)
\overline{Y}\left(\overline{Y}(e^{z_{1}D}u,x_{0})(e^{z_{2}D}v),x_{2}\right)(e^{z_{0}D}w).
\end{eqnarray*}
This implies the Jacobi identity that we need. This proves that
$(\overline{V},\overline{Y})$ carries the structure of a vertex
Leibniz algebra. Combining the $D$-bracket-derivative property and
the $D$-derivative property we see that $D$ is a derivation of
$\overline{V}$.

It can be readily seen that the map $\overline{Y}$ extends the map
$Y$, so the last assertion follows.
\end{proof}

For a vertex Leibniz algebra $V$ with $D$, we define a {\em $V$-module}
to be a module $(W,Y_{W})$ for $V$ viewed as a vertex Leibniz algebra such that (\ref{ed-module}) holds.

The proof of Proposition \ref{pleibniz-vano} (with $z_{0}=0$) immediately gives:

\bp{pbarmodule} Let $V$ be a vertex Leibniz algebra and let $(W,Y_{W})$ be a $V$-module.
Define a linear map $\overline{Y}_{W}(\cdot,x): \overline{V}\rightarrow (\End W)[[x,x^{-1}]]$  by
$$\overline{Y}_{W}(D^{n}v,x)=\left(\frac{d}{dx}\right)^{n}Y_{W}(v,x)\ \ \ \mbox{ for }n\in \N,\ v\in V.$$
Then $(W,\overline{Y}_{W})$ carries the structure of a $\overline{V}$-module.
\ep

Next, we associate a vertex algebra without vacuum to each vertex Leibniz algebra with $D$.

\bp{pquotient-skew-symmetry} Let $V$ be a vertex Leibniz algebra with $D$.
Let $J$ be the subspace of $V$, linearly spanned by vectors
\begin{eqnarray}
u_{n}v-\sum_{i\ge 0}(-1)^{n+i-1}\frac{1}{i!}D^{i}(v_{n+i}u)
\end{eqnarray}
for $u,v\in V,\ n\in \Z$.  Then $J$ is a two-sided ideal and
\begin{eqnarray}
Y(a,x)v=0\ \ \ \mbox{for }a\in J,\ v\in V.
\end{eqnarray}
Furthermore, $V/J$ is a vertex algebra without vacuum.\ep

\begin{proof} Note that for $u,v\in V,\ n\in \Z$, we have
$$\Res_{x}x^{n}\left(Y(u,x)v-e^{xD}Y(v,-x)u\right)
=u_{n}v-\sum_{i\ge 0}(-1)^{n+i-1}\frac{1}{i!}D^{i}(v_{n+i}u).$$ Let
$a,u,v\in V$. Using the $D$-bracket-derivative property and
commutator formula, we get
\begin{eqnarray*}
&&Y(a,x_{1})(Y(u,x_{2})v-e^{x_{2}D}Y(v,-x_{2})u)\\
&=&Y(a,x_{1})Y(u,x_{2})v-e^{x_{2}D}Y(a,x_{1}-x_{2})Y(v,-x_{2})u\\
&=&Y(u,x_{2})Y(a,x_{1})v+
\Res_{x_{0}}x_{1}^{-1}\delta\left(\frac{x_{2}+x_{0}}{x_{1}}\right)Y(Y(a,x_{0})u,x_{2})v\\
&&-e^{x_{2}D}Y(a,x_{1}-x_{2})Y(v,-x_{2})u.
\end{eqnarray*}
Furthermore, from the Jacobi identity for the triple $(a,v,u)$
we get
\begin{eqnarray*}
&&Y(a,x_{1}-x_{2})Y(v,-x_{2})u\nonumber\\
&=&Y(Y(a,x_{1})v,-x_{2})
+\Res_{x_{0}}x_{1}^{-1}\delta\left(\frac{x_{2}+x_{0}}{x_{1}}\right)Y(v,-x_{2})Y(a,x_{0})u.
\end{eqnarray*}
Then we obtain
\begin{eqnarray*}
&&Y(a,x_{1})(Y(u,x_{2})v-e^{x_{2}D}Y(v,-x_{2})u)\\
&=&Y(u,x_{2})Y(a,x_{1})v-e^{x_{2}D}Y(Y(a,x_{1})v,-x_{2})\\
&&+
\Res_{x_{0}}x_{2}^{-1}\delta\left(\frac{x_{1}-x_{0}}{x_{2}}\right)
\left(Y(Y(a,x_{0})u,x_{2})v
-e^{x_{2}D}Y(v,-x_{2})Y(a,x_{0})u\right).
\end{eqnarray*}
It follows that $J$ is a left ideal of $V$. Furthermore, for $a\in
J,\ v\in V$ we have
\begin{eqnarray*}
Y(a,x)v=\left(Y(a,x)v-e^{xD}Y(v,-x)a\right)+e^{xD}Y(v,-x)a.
\end{eqnarray*}
{}From this it follows that $J$ is also a right ideal. On the other
hand, by the $D$-derivative property and (\ref{epre-skew-symmetry})
we have
\begin{eqnarray*}
&&Y\left(Y(u,x_{0})v-e^{x_{0}D}Y(v,-x_{0})u,x_{2}\right)a\\
&=&Y(Y(u,x_{0})v,x_{2})a-Y(Y(v,-x_{0})u,x_{2}+x_{0})a\\
&=&0.
\end{eqnarray*}
This proves the first assertion. Then the quotient $V/J$ is a vertex
Leibniz algebra equipped with a linear operator induced from $D$,
satisfying the $D$-bracket-derivative property and the skew
symmetry. Thus $V/J$ is a vertex algebra without vacuum.
\end{proof}

Using the proof of Proposition \ref{pvla-vawv} we immediately have:

\bp{pquotient-module} Let $V$ be a vertex Leibniz algebra with $D$.
Then every $V$-module is naturally a $V/J$-module.
\ep

Let $V$ be a vertex Leibniz algebra. {}From Proposition
\ref{pleibniz-vano}, we have a vertex Leibniz algebra
$\overline{V}=\C[D]\otimes V$ with $D$. Let $\overline{J}$ denote the (two-sided)
ideal of $\overline{V}$, introduced in Proposition
\ref{pquotient-skew-symmetry}. Furthermore, by Proposition
\ref{pquotient-skew-symmetry}  the quotient algebra
$\overline{V}/\overline{J}$ is a vertex algebra without vacuum.

\bd{dextension-leibniz} {\em Let $V$ be a vertex Leibniz algebra. We
define $\widetilde{V}$ to be the vertex-algebra-without-vacuum
$\overline{V}/\overline{J}$.} \ed

\br{rmoduleisomorphism} {\em Let $V$ be a vertex Leibniz algebra.
{}From Proposition \ref{pbarmodule},  the category of $V$-modules is canonically isomorphic
to the category of $\overline{V}$-modules. On the other hand,
{}from Proposition \ref{pquotient-module} the category of $\overline{V}$-modules is
canonically isomorphic to the category of $\widetilde{V}$-modules.}
\er

As $V$ is a vertex Leibniz subalgebra of $\overline{V}$ and
$\widetilde{V}$ is a quotient algebra of $\overline{V}$, we have a
canonical homomorphism of vertex Leibniz algebras
\begin{eqnarray}
\pi_{V}: V\rightarrow \widetilde{V}.
\end{eqnarray}
The following is a universal property of $\widetilde{V}$:

\bp{pembedding} Let $V$ be a vertex Leibniz algebra as before, let
$K$ be any vertex algebra without vacuum and let $\psi: V\rightarrow
K$ be any homomorphism of vertex Leibniz algebras. Then there exists
a unique homomorphism $\tilde{\psi}$ of vertex algebras without
vacuum from $\widetilde{V}$ to $K$ such that $\tilde{\psi}\circ
\pi_{V}=\psi$. \ep

\begin{proof} It is clear that $\widetilde{V}$ as a vertex algebra without vacuum is generated
by $\pi_{V}(V)$, so the uniqueness follows immediately. We now
establish the existence. First, we define a linear map $\bar{\psi}:
\overline{V}\rightarrow K$ by $\bar{\psi}(D^{n}v)=D^{n}\psi(v)$ for
$n\in \N,\ v\in V$, where we abuse $D$ also for the $D$-operator of
$K$. This is a homomorphism of vertex Leibniz algebras, commuting
with $D$. Then $\bar{\psi}$ reduces to a homomorphism of vertex
algebras without vacuum from $\widetilde{V}$ to $K$.
\end{proof}

As an immediate consequence we have (cf. Proposition
\ref{pfaithful}):

\bc{cpiV-injective} Let $V$ be a vertex Leibniz algebra. Then homomorphism $\pi_{V}$ is injective if and
only if $V$ is a vertex Leibniz subalgebra of some vertex algebra
and if and only if $V$ admits a faithful representation. \ec

We also have:

\bp{ppiV-JV}
Let $V$ be a vertex Leibniz algebra. Then $J_{V}=\ker \pi_{V}$.
\ep

\begin{proof} Noticing that homomorphism $\pi_{V}$ from $V$ to $\widetilde{V}$ gives rise to a representation of $V$ on $\widetilde{V}$,
{}from definition we have $J_{V}\subset \ker \pi_{V}$. On the other hand, from Proposition \ref{pJV-I} $V/J_{V}$ can be embedded into
a vertex algebra $K$, so we have a homomorphism of vertex Leibniz algebras $f: V\rightarrow K$ with $\ker f=J_{V}$.
By Proposition \ref{pembedding}, there exists a homomorphism of vertex algebras without vacuum
$\psi: \widetilde{V}\rightarrow K$ such that
$f=\psi\circ \pi_{V}$.  Then $\ker \pi_{V}\subset \ker f=J_{V}$. Therefore, $\ker\pi_{V}= J_{V}$.
\end{proof}

We end this section with the following observation:

\bc{cnecessary} Let $V$ be a vertex Leibniz algebra.
A necessary condition that $V$ can be embedded into a vertex algebra is that for $u,v\in V$,
$u_{m}v=(-1)^{m-1}v_{m}u$,  where $m\in \Z$ such that $v_{m+j}u=0$ for all $j>0$.
\ec

\begin{proof}
This is because
$$u_{m}v+(-1)^{m}v_{m}u=u_{m}v-(-1)^{m-1}v_{m}u-\sum_{j\ge 1}(-1)^{m+j-1}\frac{1}{j!}D^{j}v_{m+j}u,$$
which lies in $\ker \pi_{V}\ (=J_{V})$.
\end{proof}

\section{Construction of vertex Leibniz algebras}

In this section, we first establish an analog of a theorem of Xiaoping Xu
and then using this we associate a vertex Leibniz algebra $V_{\g}$
to each Leibniz algebra $\g$ and associate a $V_{\g}$-module to
every $\g$-module.

The following is a vertex Leibniz algebra analog of a theorem of Xu (see \cite{xu}; cf.
\cite{ll}):

\bt{tanalogue} Let $V$ be a vector space equipped with a linear map
\begin{eqnarray*}
Y(\cdot,x):&& V\rightarrow \Hom (V,V((x)))\subset (\End V)[[x,x^{-1}]],\\
&&v\mapsto Y(v,x)=\sum_{n\in \Z} v_{n}x^{-n-1}
\end{eqnarray*}
and let $U$ be a subset of $V$. Assume that all the following
conditions are satisfied: For $u,v\in U$, there exists a nonnegative
integer $k$ such that
$$(x_{1}-x_{2})^{k}Y(u,x_{1})Y(v,x_{2})=(x_{1}-x_{2})^{k}Y(v,x_{2})Y(u,x_{1});$$
$V$ is linearly spanned by vectors
$$u^{(1)}_{m_{1}}\cdots u^{(r)}_{m_{r}}u$$
for $u^{(1)},\dots,u^{(r)},\; u\in U,\;  m_{1},\dots,m_{r}\in \Z$
with $r\in \N$; and for all $u\in U,\; v\in V$,
\begin{eqnarray}\label{eiterate-thm}
& &Y(Y(u,x_{0})v,x_{2})\nonumber\\
&=&\Res_{x_{1}}\left(x_{0}^{-1}\delta\left(\frac{x_{1}-x_{2}}{x_{0}}\right)
Y(u,x_{1})Y(v,x_{2})-
x_{0}^{-1}\delta\left(\frac{x_{2}-x_{1}}{-x_{0}}\right)
Y(v,x_{2})Y(u,x_{1})\right).\nonumber\\
\end{eqnarray}
Then $(V,Y)$ carries the structure of a vertex Leibniz algebra. \et

\begin{proof} Set
$$\tilde{U}=\{ Y(u,x)\;|\; u\in U\}\subset \Hom (V,V((x))).$$
{}From assumption, $\tilde{U}$ is a local subset. Then by a result of
\cite{li-local}, $\tilde{U}$ generates a vertex algebra
$\<\tilde{U}\>$, where $1_{V}$ serves as the vacuum vector and the
vertex operator map, denoted by $Y_{\E}(\cdot,x)$, is given by
\begin{eqnarray}\label{eaction-formula}
&&Y_{\E}(\alpha(x),x_{0})\beta(x)\nonumber\\
&=&\Res_{x_{1}}\left(x_{0}^{-1}\delta\left(\frac{x_{1}-x}{x_{0}}\right)
\alpha(x_{1})\beta(x)-
x_{0}^{-1}\delta\left(\frac{x-x_{1}}{-x_{0}}\right)
\beta(x)\alpha(x_{1})\right)
\end{eqnarray}
for $\alpha(x),\beta(x)\in \Hom (V,V((x)))$. It was also proved therein
that $\<\tilde{U}\>$ is also a local subspace. For $a\in U,\; v\in
V$, by (\ref{eiterate-thm}) and (\ref{eaction-formula}) we have
\begin{eqnarray}\label{eu-property}
&&Y(Y(a,x_{0})v,x)\nonumber\\
&=&\Res_{x_{1}}\left(x_{0}^{-1}\delta\left(\frac{x_{1}-x}{x_{0}}\right)
Y(a,x_{1})Y(v,x)-
x_{0}^{-1}\delta\left(\frac{x-x_{1}}{-x_{0}}\right)
Y(v,x)Y(a,x_{1})\right)\nonumber\\
 &=&Y_{\E}(Y(a,x), x_{0})Y(v,x).
 \end{eqnarray}
It follows from this and the span assumption that $Y(\cdot,x)$ maps $V$ into
$\<\tilde{U}\>$. Since $\<\tilde{U}\>$ is a local subspace of $\Hom
(V,V((x)))$, $\{ Y(v,x)\;|\; v\in V\}$ as a subset of
$\<\tilde{U}\>$ is local. Recall that Jacobi identity amounts to weak commutativity and weak associativity.
Now it remains to establish weak
associativity, which will be achieved in the following by
induction.

Let $K$ consist of each $u\in V$ such that (\ref{eiterate-thm})
holds for every $v\in V$. Now, let $a,b\in K,\ v\in V$. From
(\ref{eiterate-thm}) with $(a,b)$ in the places of $(u,v)$, there
exists a nonnegative integer $l$ such that
\begin{eqnarray*}
(x_{0}+y_{0})^{l}Y(Y(a,x_{0})b,y_{0})v
=(x_{0}+y_{0})^{l}Y(a,x_{0}+y_{0})Y(b,y_{0})v.
\end{eqnarray*}
Using this, (\ref{eu-property}), and the weak associativity for
$\<\tilde{U}\>$,  replacing $l$ with a large one if necessary, we
have
\begin{eqnarray*}
&&(x_{0}+y_{0})^{l}Y(Y(Y(a,x_{0})b,y_{0})v,z_{0})\\
&=&(x_{0}+y_{0})^{l}Y\left(Y(a,x_{0}+y_{0})Y(b,y_{0})v,z_{0}\right)\\
&=&(x_{0}+y_{0})^{l}Y_{\E}\left(\overline{Y(a,z_{0})},x_{0}+y_{0}\right)
\overline{Y\left(Y(b,y_{0})v,z_{0}\right)}\\
&=&(x_{0}+y_{0})^{l}Y_{\E}\left(\overline{Y(a,z_{0})},x_{0}+y_{0}\right)
Y_{\E}\left(\overline{Y(b,z_{0})},y_{0}\right)\overline{Y(v,z_{0})}\\
&=&(x_{0}+y_{0})^{l}Y_{\E}\left(Y_{\E}\left(\overline{Y(a,z_{0})},x_{0}\right)
\overline{Y(b,z_{0})},y_{0}\right)\overline{Y(v,z_{0})}\\
&=&(x_{0}+y_{0})^{l}
Y_{\E}\left(\overline{Y\left(Y(a,x_{0})b,z_{0}\right)},y_{0}\right)
\overline{Y(v,z_{0})},
 \end{eqnarray*}
where $\overline{X}=X$ for all the bar objects; the only purpose is
to make the equation easier to read. Multiplying both sides by
formal series $(y_{0}+x_{0})^{-l}$ we obtain
\begin{eqnarray*}
Y(Y(Y(a,x_{0})b,y_{0})v,z_{0})
=Y_{\E}\left(\overline{Y\left(Y(a,x_{0})b,z_{0}\right)},y_{0}\right)
\overline{Y(v,z_{0})}.
 \end{eqnarray*}
Furthermore, we have
\begin{eqnarray*}
&&Y(Y(Y(a,x_{0})b,y_{0})v,z_{0})\\
&=&Y_{\E}\left(\overline{Y(Y(a,x_{0})b,z_{0})},y_{0}\right)
\overline{Y(v,z_{0})}\\
&=&\Res_{x_{1}}y_{0}^{-1}\delta\left(\frac{x_{1}-z_{0}}{y_{0}}\right)
Y(Y(a,x_{0})b,x_{1})Y(v,z_{0})\\
&&-\Res_{x_{1}}y_{0}^{-1}\delta\left(\frac{z_{0}-x_{1}}{-y_{0}}\right)
Y(v,z_{0})Y(Y(a,x_{0})b,x_{1}).
\end{eqnarray*}
This shows that $a_{m}b\in K$ for $m\in \Z$. It follows from
induction and the span assumption that $V=K$. This proves that
(\ref{eiterate-thm}) holds for {\em all } $u,v\in V$. Then weak
associativity follows. Therefore, $(V,Y)$ carries the structure of a
vertex Leibniz algebra.
\end{proof}

\br{rclassical-notion} {\em  Recall that a (left) {\em
Leibniz algebra} is a non-associative algebra $A$ satisfying the
Leibniz role:
$$a(bc)=(ab)c+b(ac)\ \ \mbox{ for }a,b,c\in A.$$
Let $J$ be the subspace linearly spanned by $a^{2}$ for $a\in A$. A
simple fact is that $J\cdot A=0$ and $A\cdot J\subset J$, which
implies that $J$ is a two-sided ideal of $A$. Furthermore, the
quotient algebra $A/J$ is a Lie algebra. } \er

Now, let $\g$ be a Leibniz algebra. Form the Loop algebra
$L(\g)=\g\otimes \C[t,t^{-1}]$, where for $a,b\in \g,\ m,n\in \Z$,
\begin{eqnarray}
[a\otimes t^{m},b\otimes t^{n}]=[a,b]\otimes t^{m+n}.
\end{eqnarray}
It can be readily seen that the loop algebra
 $L(\g)$ is also a Leibniz algebra.

\br{raffinization-leibniz} {\em Let $\g$ be a Leibniz algebra and
let $\<\cdot,\cdot\>$ be a bilinear form on $\g$. Consider the
affinization $\hat{\g}=\g[t,t^{-1}]\oplus \C {\bf k}$, where
\begin{eqnarray*}
&&[{\bf k}, \hat{\g}]=0=[\hat{\g},{\bf k}],\\
&&[a\otimes t^{m},b\otimes t^{n}]=[a,b]\otimes
t^{m+n}+m\<a,b\>\delta_{m+n,0}{\bf k}
\end{eqnarray*}
for $a,b\in \g,\ m,n\in \Z$. It is straightforward to show that
$\hat{\g}$ is a Leibniz algebra if and only if
\begin{eqnarray}
\<a,[b,c]\>=\<[a,b],c\>=-\<b,[a,c]\>\ \ \ \mbox{ for }a,b,c\in \g.
\end{eqnarray}
{}From this, we see that if $\<\cdot,\cdot\>$ is non-degenerate,
then $\hat{\g}$ is a Leibniz algebra if and only if $\g$ is a Lie
algebra and $\<\cdot,\cdot\>$ is an invariant bilinear form.} \er

Denote by $L(\g)_{Lie}$ the associated Lie algebra. Recall that
vectors $[u,u]$ for $u\in L(\g)$ linearly span a two-sided ideal of
$L(\g)$ and  $L(\g)_{Lie}$ is the quotient algebra.

\br{rtwo-same}
{\em Note that we have two Lie algebras $L(\g)_{Lie}$ and $L(\g_{Lie})$.
We here show that the two Lie algebras are actually isomorphic.
Let $J(\g)$ and $J(L(\g))$ denote the ideals of $\g$ and $L(\g)$, respectively.
We need to show that $J(L(\g))=J(\g)\otimes \C[t,t^{-1}].$
For any $a,b\in \g,\; m,n\in \Z$, we have
\begin{eqnarray*}
[a\otimes t^{m},b\otimes t^{n}]+[b\otimes t^{n},a\otimes t^{m}]
=([a,b]+[b,a])\otimes t^{m+n}
\in J(\g)\otimes \C[t,t^{-1}].
\end{eqnarray*}
It follows that $J(L(\g))\subset J(\g)\otimes \C[t,t^{-1}]$.
On the other hand, for $a\in \g,\; m\in \Z$, we have
$$[a,a]\otimes t^{m}=\frac{1}{2}\left([a\otimes t^{m}, a\otimes 1]+[a\otimes 1,a\otimes t^{m}]\right)\in J(L(\g)).$$
This implies $J(\g)\otimes \C[t,t^{-1}]\subset J(L(\g))$. Thus
$J(L(\g))=J(\g)\otimes \C[t,t^{-1}].$ } \er

A {\em left $\g$-module} is a vector space $W$ on which $\g$ acts
from left such that
\begin{eqnarray}
a(bw)-b(aw)=[a,b]w\ \ \ \mbox{ for }a,b\in \g,\; w\in W.
\end{eqnarray}
It can be readily seen that a left $\g$-module exactly amounts to a module for the associated Lie algebra $\g_{Lie}$.
Note that a Leibniz algebra $\g$ itself is always a left $\g$-module.

Let $U$ be a (left) $\g$-module, which amounts to a $\g_{Lie}$-module.
Let $\g_{Lie}\otimes t\C[t]$ act trivially on $U$, to make $U$ an $L(\g_{Lie})^{+}$-module, where
$$L(\g_{Lie})^{+}=\g_{Lie}\otimes \C[t]=\g_{Lie}+ (\g_{Lie}\otimes t\C[t]).$$
Form an induced module
\begin{eqnarray}
V_{\g}(U)=U(L(\g_{Lie}))\otimes _{U(L(\g_{Lie})^{+})}U.
\end{eqnarray}
It follows that $V_{\g}(U)$ is a {\em restricted} $L(\g)$-module in the sense that for any $a\in \g,\ w\in V_{\g}(U)$,
$(a\otimes t^{n})w=0$ for $n$ sufficiently large. Set
\begin{eqnarray}
V_{\g}=V_{\g}(\g)=U(L(\g_{Lie}))\otimes _{U(L(\g_{Lie})^{+})}\g.
\end{eqnarray}
 We have:

\bt{tloop-leibniz}
 There exists a vertex Leibniz algebra structure on $V_{\g}$, which is uniquely
determined by the condition that $Y(a,x)=a(x)$ for $a\in \g$, where
$$a(x)=\sum_{n\in \Z}(a\otimes t^{n})x^{-n-1}.$$
Furthermore, for any restricted $L(\g)$-module $W$, there is
a $V_{\g}$-module structure $Y_{W}(\cdot,x)$ on $W$, which is uniquely determined by
the condition that
$$Y_{W}(a,x)=a(x)\ \ \ \mbox{ for }a\in \g.$$
\et

\begin{proof} Let $W$ be any restricted $L(\g)$-module. Set
$$U_{W}=\{ a(x)\;|\; a\in \g\}\subset \E(W)\ \left(=\Hom (W,W((x)))\right).$$
For any $a,b\in \g$, from the commutation relation of $L(\g)$ we have
\begin{eqnarray}\label{eab-comm}
[a(x_{1}),b(x_{2})]=[a,b](x_{2})x_{1}^{-1}\delta\left(\frac{x_{2}}{x_{1}}\right),
\end{eqnarray}
which implies
$$(x_{1}-x_{2})[a(x_{1}),b(x_{2})]=0.$$
Thus $U_{W}$ is a local subset of $\E(W)$.
By a result of \cite{li-local}, $U_{W}$ generates a vertex algebra $\<U_{W}\>$ where
the vertex operator map is denoted by $Y_{\E}(\cdot,x)$ and the identity operator $1_{W}$ is
 the vacuum vector. Furthermore, $W$ is a faithful $\<U_{W}\>$-module
 with $Y_{W}(u(x),z)=u(z)$ for $u(x)\in \<U_{W}\>$.
With (\ref{eab-comm}), from \cite{li-local} (Lemma 2.3.5) we have
\begin{eqnarray}
[Y_{\E}(a(x),x_{1}),Y_{\E}(b(x),x_{2})]=Y_{\E}([a,b](x),x_{2})x_{1}^{-1}\delta\left(\frac{x_{2}}{x_{1}}\right).
\end{eqnarray}
This shows that $\<U_{W}\>$ is an $L(\g)$-module with $a(z)=Y_{\E}(a(x),z)$ for $a\in \g$. Furthermore,
$$a(x)_{0}b(x)=[a,b](x),\ \ \ a(x)_{n}b(x)=0\ \ \ \mbox{ for }n\ge 1. $$
Then it follows from the construction of $V_{\g}$ that there exists an $L(\g)$-module homomorphism
$\psi_{W}$ from $V_{\g}$ to $\<U_{W}\>$, sending $a$ to $a(x)$ for $a\in \g$.

Now, we specialize $W=V_{\g}$ and we denote the corresponding $L(\g)$-module homomorphism
$\psi_{W}$ by $\psi_{x}$.
For $v\in V_{\g}$, set
$$Y(v,x)=\psi_{x}(v)\in \<U_{V_{\g}}\>\subset \Hom (V,V((x))).$$
We have
$$Y(a,x)=\psi_{x}(a)=a(x)\ \ \ \mbox{ for }a\in \g$$
and furthermore, for $v\in V_{\g}$,
\begin{eqnarray*}
&&Y(Y(a,x_{0})v,x_{2})=Y(a(x_{0})v,x_{2})
=\psi_{x_{2}}(a(x_{0})v)\\
&=&Y_{\E}(a(x),x_{0})\psi_{x_{2}}(v)\\
&=&\Res_{x_{1}}\left(x_{0}^{-1}\delta\left(\frac{x_{1}-x_{2}}{x_{0}}\right)
a(x_{1})\psi_{x_{2}}(v)-
x_{0}^{-1}\delta\left(\frac{x_{2}-x_{1}}{-x_{0}}\right)
\psi_{x_{2}}(v)a(x_{1})\right)\\
&=&\Res_{x_{1}}\left(x_{0}^{-1}\delta\left(\frac{x_{1}-x_{2}}{x_{0}}\right)
Y(a,x_{1})Y(v,x_{2})-
x_{0}^{-1}\delta\left(\frac{x_{2}-x_{1}}{-x_{0}}\right)
Y(a,x_{2})Y(u,x_{1})\right).
\end{eqnarray*}
Now the first assertion immediately follows from Theorem \ref{tanalogue} where the uniqueness
follows from the spanning property.

Next, we come back to a general restricted $L(\g)$-module $W$.
Recall that $\psi_{W}$ is an $L(\g)$-module homomorphism from $V_{\g}$ to $\<U_{W}\>$
with $\psi_{W}(a)=a(x)$ for $a\in \g$. It follows that $\psi_{W}$
is a homomorphism of vertex Leibniz algebras.
As $W$ is a $\<U_{W}\>$-module, $W$ becomes a $V_{\g}$-module through $\psi_{W}$.
\end{proof}

\br{rspecial-case} {\em Let $\g$ be a Lie algebra, which is a
Leibniz algebra with $\g_{Lie}=\g$. It follows from the
constructions of $V_{\g}$ and $V_{\hat{\g}}(0,0)$ (recall Example \ref{eaxmple-verma}) and the P-B-W
theorem that $V_{\g}$ is isomorphic to the $L(\g)$-submodule of
$V_{\hat{\g}}(0,0)$, generated by $\g$. As $Y(a,x)=a(x)$ for $a\in
\g$, it follows that $V_{\g}$ is a closed subspace. Then $V_{\g}$ is
a vertex Leibniz subalgebra. Furthermore, for $a,b\in \g\subset
V_{\hat{\g}}(0,0)$, we have
$$[a,b](-2){\bf 1}=[a(-1),b(-1)]{\bf 1}=a(-1)b-b(-1)a\in V_{\g}.$$
Recall that $\D u=u(-2){\bf 1}$ for $u\in \g$. We see that if
$[\g,\g]=\g$, then $\D V_{\g}\subset V_{\g}$ and $V_{\g}$ is a
vertex algebra without vacuum.} \er

\bp{psubva-condition} Let $\g$ be a Leibniz algebra and let $V_{\g}$
be the associated vertex Leibniz algebra. Then $V_{\g}$ is a vertex
Leibniz subalgebra of an ordinary vertex algebra if and only if $\g$
is a Lie algebra. \ep

\begin{proof} If $\g$ is a Lie algebra, from Remark \ref{rspecial-case}
we see that $V_{\g}$ is isomorphic to the
vertex Leibniz subalgebra of vertex algebra $V_{\hat{\g}}(0,0)$, generated
by $\g$.

Now, assume that there exists a vertex algebra $V$ containing
$V_{\g}$ as a vertex Leibniz subalgebra. For $a,b\in \g\subset
V_{\g}$, with $Y(a,x)=a(x)$ we have
\begin{eqnarray}\label{ea0banb}
a_{0}b=[a,b]\ \ \mbox{ and }\ \ a_{n}b=0\ \
\mbox{ for }n\ge 1.
\end{eqnarray}
On the other hand, as $V_{\g}$ is a vertex Leibniz subalgebra of
vertex algebra $V$, we have
$$Y(a,x)b=e^{x\D}Y(b,-x)a,$$
where $\D$ is the linear operator of $V$, defined by $\D
(v)=v_{-2}{\bf 1}$ for $v\in V$. Then
$$[a,b]=a_{0}b=\sum_{i\ge 0}(-1)^{i+1}\frac{1}{i!}D^{i}b_{i}a=-b_{0}a=-[b,a].$$
Thus $\g$ must be a Lie algebra. Therefore, if $\g$ is an authentic
Leibniz algebra (not a Lie algebra), $V_{\g}$ is not a subalgebra of
a vertex algebra.
\end{proof}

\end{document}